\documentclass[12pt]{article}
\usepackage{latexsym,amsfonts,amssymb}
\usepackage[dvips]{color}
\usepackage{colordvi,multicol}
\usepackage{amsmath}
\usepackage{graphicx}
\usepackage{pdflscape}
\usepackage{rotating}

\usepackage{latexsym,amsfonts,amssymb}
\usepackage[dvips]{color}
\usepackage{colordvi,multicol}
\usepackage{amsmath}
\usepackage{graphicx}
\usepackage{pdflscape}
\usepackage{rotating}
\usepackage{breqn}

\def \cal{\mathcal}

\newtheorem{thm}{Theorem}[section]

\newtheorem{pro}[thm]{Proposition}

\newtheorem{exa}[thm]{Example}

\date{}
\begin{document}

\title{\bf On the measure concentration of infinitely divisible distributions}
\author{Jing Zhang, Ze-Chun Hu and Wei Sun\thanks{Corresponding author.}\\ \\
  {\small School of Mathematics and Statistics, Hainan Normal University, Haikou 571158, China}\\ \\
 {\small College of Mathematics, Sichuan University, Chengdu 610065, China}\\ \\
{\small Department of Mathematics and Statistics, Concordia University, Montreal H3G 1M8,  Canada}\\ \\
{\small  zh\_jing0820@hotmail.com\ \ \ \ zchu@scu.edu.cn\ \ \ \ wei.sun@concordia.ca}}

\maketitle

\begin{abstract}
\noindent Let ${\cal I}$ be the set of all infinitely divisible random variables\ with finite second moments, ${\cal I}_0=\{X\in{\cal I}:{\rm Var}(X)>0\}$, $P_{\cal I}=\inf_{X\in{\cal I}}P\{|X-E[X]|\le \sqrt{{\rm Var}(X)}\}$ and $P_{{\cal I}_0}=\inf_{X\in{\cal I}_0} P\{|X-E[X]|< \sqrt{{\rm Var}(X)}\}$. Firstly,  we prove that $P_{{\cal I}}\ge P_{{\cal I}_0}>0$. Secondly, we find the exact values of $\inf_{X\in{\cal J}}P\{|X-E[X]|\le \sqrt{{\rm Var}(X)}\}$ and $\inf_{X\in\cal J} P\{|X-E[X]|< \sqrt{{\rm Var}(X)}\}$ for the cases that $\cal J$ is the set of all geometric random variables, symmetric geometric random variables, Poisson random variables and symmetric Poisson random variables, respectively.  As a consequence, we obtain that $P_{\cal I}\le e^{-1}\sum_{k=0}^{\infty}\frac{1}{2^{2k}(k!)^2}\approx 0.46576$ and $P_{{\cal I}_0}\le e^{-1}\approx 0.36788$.
\end{abstract}



\section{Introduction}

A distribution $\mu$ on $\mathbb{R}$ is infinitely divisible
if it can be expressed as the distribution of the sum of an arbitrary number of i.i.d. random variables. It is well known that each L\'evy process can be associated with an infinitely divisible distribution. Infinitely divisible distributions play a fundamental role in probability theory and stochastic processes. They have found applications in various fields, including physics, chemistry, climate changes, communications  and finance. Although the study of infinitely divisible distributions  has a long history, there are still many important related problems remain unsolved, e.g., Getoor's conjecture  that essentially all L\'evy processes satisfy Hunt's hypothesis (H) (cf. \cite{H1, H}). It is worth mentioning that there is a close connection between infinite divisibility and some challenging problems in other math fields, e.g., the recently discovered equivalence between the Riemann hypothesis and infinite divisibility (cf. \cite{N1, N}).

In \cite{S1,S2,S3}, we initiate the study of the variation comparison between infinitely divisible distributions and the normal distribution. Let $X$ be a random variable with finite second moment. We consider the inequality:
\begin{equation}\label{27a}
P\left\{|X-E[X]|\le \sqrt{{\rm Var}(X)}\right\}\ge P\{|Z|\le 1\},
\end{equation}
 where $Z$ is a standard normal random variable. We prove that this inequality holds for many familiar infinitely divisible continuous distributions including the Gamma, Laplace, Gumbel, Logistic, Pareto, infinitely divisible Weibull,  log-normal, student's $t$, inverse Gaussian and $F$-distributions. In \cite{S2}, we also discuss the quantity $P\{|X-E[X]|\le \sqrt{{\rm Var}(X)}\}$ for some infinitely divisible discrete distributions. Numerical results show that the variation comparison inequality (\ref{27a}) does not hold for some negative binomial distributions or Poisson distributions.

In this paper, we will further discuss the variation $P\{|X-E[X]|\le \sqrt{{\rm Var}(X)}\}$. Define
$$
{\cal D}:=\{X:\, X\ {\rm is\ a\ random\ variable\ with}\ E[X^2]<\infty\}.
$$
First, we point out that
$$
\inf\limits_{X\in{\cal D}} P\left\{|X-E[X]|\le \sqrt{{\rm Var}(X)}\right\}=0.
$$
In fact, this can be seen from the following simple example.
\begin{exa}\label{exa1}
Fix $a_1,a_2\in\mathbb{R}$ with $0\le a_1<a_2$ and let $\varepsilon\in(0,1]$ be arbitrary. Suppose that $X$ is a discrete random variable with the probability mass function:
$$
P\{X=a_1\}=P\{X=-a_1\}=\frac{\varepsilon}{2},\ \ \ \ P\{X=a_2\}=P\{X=-a_2\}=\frac{1-\varepsilon}{2}.
$$
We have that $E[X]=0$ and
$$
{\rm Var}(X)=a_1^2\varepsilon+a_2^2(1-\varepsilon)\in [a_1^2,a_2^2).
$$
Then,
$$
P\left\{|X-E[X]|\le \sqrt{{\rm Var}(X)}\right\}=\varepsilon.
$$
Since $\varepsilon$ is arbitrary, the variation $
P\{|X-E[X]|\le \sqrt{{\rm Var}(X)}\}$ can take any value in $(0,1]$.
\end{exa}

Different from general distributions, we discover that infinitely divisible distributions exhibit an interesting measure concentration phenomenon. Define
$$
{\cal I}:=\{X:\, X\ {\rm is\ an\ infinitely\ divisible\ random\ variable\ with}\ E[X^2]<\infty\},
$$
$$
{\cal I}_0:=\{X\in {\cal I} :\, {\rm Var}(X)>0\},
$$
and
$$
P_{\cal I}:=\inf\limits_{X\in{\cal I}} P\left\{|X-E[X]|\le \sqrt{{\rm Var}(X)}\right\},\ \ \ \ P_{{\cal I}_0}:=\inf\limits_{X\in{\cal I}_0} P\left\{|X-E[X]|< \sqrt{{\rm Var}(X)}\right\}.
$$
In the next section,  we will show that $P_{\cal I}\ge P_{{\cal I}_0}>0$ (see Theorem \ref{thm1}).

In Section 3, we determine the values of  $\inf_{X\in \cal J} P\{|X-E[X]|\le \sqrt{{\rm Var}(X)}\}$ and $\inf_{X\in \cal J}P\{|X-E[X]|< \sqrt{{\rm Var}(X)}\}$ when $\cal J$ stands for all geometric random variables or symmetric geometric random variables with parameter $p\in (0,1)$ (see Propositions \ref{pro1a} and \ref{pro1b}). In Section 4, we determine the values of $\inf_{X\in \cal J} P\{|X-E[X]|\le \sqrt{{\rm Var}(X)}\}$ and $\inf_{X\in \cal J} P\{|X-E[X]|< \sqrt{{\rm Var}(X)}\}$ when $\cal J$ stands for all Poisson random variables or symmetric Poisson random variables with parameter $\lambda>0$ (see Theorems \ref{pro1c}-\ref{thm2}). Theorem \ref{pro2}  gives an upper bound for $P_{{\cal I}_0}$ and Theorem \ref{thm2}  gives an upper bound for   $P_{\cal I}$.

So far we have not been able to determine the exact values of $P_{\cal I}$ and $P_{{\cal I}_0}$. We leave them as an open problem for the community of Probability and Statistics. Throughout this paper, we use $N_{\lambda}$ to denote a Poisson random variable with parameter $\lambda>0$.

\section{$P_{\cal I}\ge P_{{\cal I}_0}>0$}\setcounter{equation}{0}

\begin{thm}\label{thm1} We have $P_{\cal I}\ge P_{{\cal I}_0}>0$. Moreover, there exists  $Y\in P_{{\cal I}_0}$ such that $P_{{\cal I}_0}=P\{|Y-E[Y]|< \sqrt{{\rm Var}(Y)}\}$.
\end{thm}

\noindent Proof.\ \ We choose $\{X_n\in{{\cal I}_0}\}_{n=1}^{\infty}$ satisfying
\begin{eqnarray}\label{11}
P\left\{|X_n-E[X_n]|< \sqrt{{\rm Var}(X_n)}\right\}< P_{{\cal I}_0}+\frac{1}{n}.
\end{eqnarray}
Define
$$
Y_n=\frac{X_n-E[X_n]}{\sqrt{{\rm Var}(X_n)}}.
$$
Then, we have $Y_n\in {\cal I}_0$ and ${\rm Var}(Y_n)=1$. Thus,  $\{Y_n\}_{n=1}^{\infty}$ is uniformly integrable.

Denote by $\mu_n$ the distribution of $Y_n$. By ${\rm Var}(Y_n)=1$, we know that $\{\mu_n\}_{n=1}^{\infty}$ is tight. Without loss of generality, we assume that $\mu_n$ converges weakly to $\mu$  as $n\rightarrow\infty$. Let $Y$ be a random variable with distribution $\mu$. By the Skorohod representation theorem, we may assume without loss of generality that $Y_n$ converges to $Y$ a.s. as $n\rightarrow\infty$. By the uniform integrability of $\{Y_n\}_{n=1}^{\infty}$ and Fatou's lemma, we find that $E[Y]=0$ and ${\rm Var}(Y)\le 1$. Further, by \cite[Lemma 7.8, page 34]{S}, we know that $Y\in{\cal I}$.

Note that $P_{{\cal I}_0}<1$. By the weak convergence of $\{Y_n\}_{n=1}^{\infty}$ and (\ref{11}), we get
\begin{equation}\label{R1}
 P\left\{|Y|< 1\right\}\le \liminf_{n\rightarrow\infty}P\left\{|Y_n|< 1\right\}\le P_{{\cal I}_0}<1.
\end{equation}
Assume that $P_{{\cal I}_0}=0$. Then, by (\ref{R1}), we get $P\{|Y|< 1\}=0$. Thus, $P\{|Y|=1\}=1$ by $E[Y]=0$ and ${\rm Var}(Y)\le 1$. Hence $P\{Y=1\}= P\{Y=-1\}=\frac{1}{2}$, which contradicts  the infinite divisibility of $Y$. Therefore,
$P_{{\cal I}_0}>0.$

The inequality $P_{\cal I}\ge P_{{\cal I}_0}$ is obvious, since it is easy to see that
\begin{align*}
\inf\limits_{X\in{\cal I}} P\left\{|X-E[X]|\le \sqrt{{\rm Var}(X)}\right\}
&=\inf\limits_{X\in{\cal I}_0} P\left\{|X-E[X]|\le \sqrt{{\rm Var}(X)}\right\}\\
&\geq \inf\limits_{X\in{\cal I}_0} P\left\{|X-E[X]|<\sqrt{{\rm Var}(X)}\right\}.
\end{align*}

By $E[Y]=0$ and (\ref{R1}), we find that ${\rm Var}(Y)>0$ and hence $Y\in{\cal I}_0$. Further, by ${\rm Var}(Y)\le 1$ and (\ref{R1}), we get
$$
P\left\{|Y-E[Y]|< \sqrt{{\rm Var}(Y)}\right\}\le P_{{\cal I}_0}.
$$
Therefore,
$$
P_{{\cal I}_0}=P\left\{|Y-E[Y]|< \sqrt{{\rm Var}(Y)}\right\}.
$$
The proof is complete.\hfill \fbox

We would like to point out that  the one standard deviation plays an important role in the proof of Theorem \ref{thm1}. In fact, we have the following result.
\begin{pro} For any  $\varepsilon\in (0,\frac{\sqrt{2}}{2})$, we have
$$\inf\limits_{X\in{\cal I}} P\left\{|X-E[X]|\le \varepsilon\sqrt{{\rm Var}(X)}\right\}=0.$$
\end{pro}

\noindent Proof.\ \ Suppose $\varepsilon\in (0,\frac{\sqrt{2}}{2})$. Let $N_{\lambda}$  be a Poisson random variable with parameter $\lambda\in (\varepsilon^2,\frac{1}{2})$. Then, we have
$$
\lambda-\varepsilon\sqrt{\lambda}>0,\ \ \ \ \lambda+\varepsilon\sqrt{\lambda}<1,
$$
which implies that
\begin{eqnarray*}
 P\left\{|N_{\lambda}-E[N_{\lambda}]|\le \varepsilon\sqrt{{\rm Var}(N_{\lambda})}\right\}=P\left\{\lambda-\varepsilon\sqrt{\lambda}\le N_{\lambda}\le \lambda+\varepsilon\sqrt{\lambda}\right\}=0.
\end{eqnarray*}\hfill \fbox

Note that, by $\lim_{\lambda\downarrow 0}P\left\{|N_{\lambda}-E[N_{\lambda}]|< \sqrt{{\rm Var}(N_{\lambda})}\right\}=1$, we find that
$$
\sup\limits_{X\in{\cal I}_0} P\left\{|X-E[X]|< \sqrt{{\rm Var}(X)}\right\}=1.
$$

\section{Geometric distribution and symmetric geometric distribution}\setcounter{equation}{0}

\subsection{Geometric distribution}

Let $X_p$ be a geometric random variable  with parameter $p\in (0,1)$:
$$P\{X_p=k\}=p(1-p)^k,\ \ \ \ k=0,1,2,\dots.$$
We have
$$E[X_p]=\frac{1}{p}-1,\ \ \ \ {\rm Var}(X_p)=\frac{1-p}{p^2}.$$

\begin{pro}\label{pro1a} We have
\begin{equation}\label{as1}
\inf_{p\in (0,1)}P\left\{|X_p-E[X_p]|\le \sqrt{{\rm Var}(X_p)}\right\}=\inf_{p\in (0,1)}P\left\{|X_p-E[X_p]|< \sqrt{{\rm Var}(X_p)}\right\}=\frac{3}{4}.
\end{equation}
\end{pro}

\noindent Proof.\ \ Set $q=1-p$. Note that $\frac{1}{p}-1-\frac{\sqrt{1-p}}{p}< 0$ since
\begin{eqnarray*}
\frac{1}{p}-1-\frac{\sqrt{1-p}}{p}< 0&\Leftrightarrow&\frac{q-\sqrt{q}}{p}<0\\
&\Leftrightarrow& q<\sqrt{q}\\
&\Leftrightarrow&0< q<1.
\end{eqnarray*}
Then, we have that
\begin{eqnarray*}
P\left\{|X_p-E[X_p]|\le \sqrt{{\rm Var}(X_p)}\right\}&=&P\left\{\frac{1}{p}-1-\frac{\sqrt{1-p}}{p}\le X_p\le \frac{1}{p}-1+\frac{\sqrt{1-p}}{p} \right\}\nonumber\\
&=&P\left\{0\le X_p\le \frac{1}{p}-1+\frac{\sqrt{1-p}}{p}\right\}.
\end{eqnarray*}

\noindent (i) For $\frac{3}{4}<p<1$, we have $0<\frac{1}{p}-1+\frac{\sqrt{1-p}}{p}<1$ since
\begin{eqnarray*}
0<\frac{1}{p}-1+\frac{\sqrt{1-p}}{p}<1&\Leftrightarrow&1<\frac{1+\sqrt{1-p}}{p}<2\\
&\Leftrightarrow& p-1<\sqrt{1-p}< 2p-1\\
&\Leftrightarrow& 4p^2-3p>0\\
&\Leftrightarrow&p>\frac{3}{4}.
\end{eqnarray*}
Hence,
$$P\left\{0\le X_p\le \frac{1}{p}-1+\frac{\sqrt{1-p}}{p}\right\}=P\{X_p=0\}=p>\frac{3}{4}.$$

\noindent (ii) For $0<p\le\frac{3}{4}\Leftrightarrow \frac{1}{4}\le q<1$, we have
\begin{eqnarray*}
P\left\{0\le X_p\le \frac{1}{p}-1+\frac{\sqrt{1-p}}{p} \right\}&=&\sum_{k=0}^{\lfloor\frac{1}{p}-1+\frac{\sqrt{1-p}}{p}\rfloor}p(1-p)^k\\
&=&1-(1-p)^{\lfloor\frac{1}{p}-1+\frac{\sqrt{1-p}}{p}\rfloor+1}\\
&=&1-(1-p)^{\lfloor\frac{1+\sqrt{1-p}}{p}\rfloor}\\
&=&1-q^{\lfloor\frac{1}{1-\sqrt{q}}\rfloor}.
\end{eqnarray*}
Hereafter we use $\lfloor x\rfloor$ to denote the greatest integer less than or equal to $x$.

Define $$f(q):=1-q^{\lfloor\frac{1}{1-\sqrt{q}}\rfloor},\ \ \ \ \frac{1}{4}\le q<1.$$
Suppose  $\frac{1}{1-\sqrt{q}}\in [x,x+1)$ for some $x\in\mathbb{N}$, $x\ge 2$. Then,
$f(q)=1-q^x$ when $q\in[(\frac{x-1}{x})^2,(\frac{x}{x+1})^2)$. Since $f'(q)=-xq^{x-1}<0$, $f(q)$ is decreasing in $[(\frac{x-1}{x})^2,(\frac{x}{x+1})^2)$. Thus,
$$\inf\limits_{q\in [(\frac{x-1}{x})^2,(\frac{x}{x+1})^2)}f(q)=\lim\limits_{q\uparrow (\frac{x}{x+1})^2}f(q)=1-\left(\frac{x}{x+1}\right)^{2x}=1-\frac{1}{(1+1/x)^{2x}}.$$
Set $a_n:=1-\frac{1}{(1+1/n)^{2n}}$, $n\ge2$. Then, we have
$$\min\limits_{n\ge2}a_n=a_2\approx0.80247.$$
Hence, by cases (i) and (ii), we get
$$
\inf_{p\in (0,1)}P\left\{|X_p-E[X_p]|\le \sqrt{{\rm Var}(X_p)}\right\}=\frac{3}{4}.
$$
Therefore, the proof of (\ref{as1}) is complete by  noting that
$$
P\left\{\left|X_{\frac{3}{4}}-E\left[X_{\frac{3}{4}}\right]\right|< \sqrt{{\rm Var}\left(X_{\frac{3}{4}}\right)}\right\}=P\left\{\left|X_{\frac{3}{4}}-\frac{1}{3}\right|< \frac{2}{3}\right\}=P\left\{X_{\frac{3}{4}}=0\right\}=\frac{3}{4}.
$$
\hfill \fbox

\subsection{Symmetric geometric distribution}

Suppose that $X_p^{(1)}$ and $X_p^{(2)}$ are independent geometric random variables with the same parameter $p\in(0,1)$. Let $X^S_p=X_p^{(1)}-X_p^{(2)}$. Then,
$$E[X^S_p]=0,\ \ \ \  {\rm Var}(X^S_p)=\frac{2(1-p)}{p^2}.$$
Hence,
$$P\left\{|X^S_p-E[X^S_p]|\le \sqrt{{\rm Var}(X^S_p)}\right\}=P\left\{-\frac{\sqrt{2(1-p)}}{p}\le X^S_p\le \frac{\sqrt{2(1-p)}}{p} \right\}.
$$
We have
\begin{eqnarray*}
P\{X^S_p=k\}&=&P\{X_p^{(1)}-X_p^{(2)}=|k|\}\nonumber\\
&=&\sum_{l=0}^{\infty}P\{X_p^{(1)}=|k|+l\}P\{X_p^{(2)}=l\}\nonumber\\
&=&       \sum_{l=0}^{\infty}p(1-p)^{|k|+l}p(1-p)^{l}\\
&=&      \frac{(1-q)q^{|k|}}{1+q}.
\end{eqnarray*}

\begin{pro}\label{pro1b} We have
\begin{eqnarray}\label{as2}
&&\inf_{p\in (0,1)}P\left\{|X^S_p-E[X^S_p]|\le \sqrt{{\rm Var}(X^S_p)}\right\}=\inf_{p\in (0,1)}P\left\{|X^S_p-E[X^S_p]|< \sqrt{{\rm Var}(X^S_p)}\right\}\nonumber\\
&=&\frac{\sqrt{3}}{3}\approx  0.57735.
\end{eqnarray}
\end{pro}

\noindent Proof.\ \  (i) For $0<\frac{\sqrt{2(1-p)}}{p}<1$, i.e., $\sqrt{3}-1<p <1$, we have
\begin{eqnarray*}
P\left\{|X^S_p-E[X^S_p]|\le \sqrt{{\rm Var}(X^S_p)}\right\}&=&P\{X^S_p=0\}\\
&=&\frac{p}{2-p}\\
&>&\frac{\sqrt{3}}{3}.
\end{eqnarray*}

\noindent (ii) For $\frac{\sqrt{2(1-p)}}{p}\ge1$, i.e., $0<p \le\sqrt{3}-1$, we have
\begin{eqnarray*}
P\left\{|X^S_p-E[X^S_p]|\le \sqrt{{\rm Var}(X^S_p)}\right\}
&=&\frac{1-q}{1+q}+2\sum_{k=1}^{\lfloor\frac{\sqrt{2(1-p)}}{p}\rfloor}\frac{(1-q)q^k}{1+q}\\
&=&\frac{1+q-2q^{\lfloor\frac{\sqrt{2(1-p)}}{p}\rfloor+1}}{1+q}\\
&=&1-\frac{2q^{\lfloor\frac{\sqrt{2q}}{1-q}\rfloor+1}}{1+q}.
\end{eqnarray*}

Define $$f(q):=1-\frac{2q^{\lfloor\frac{\sqrt{2q}}{1-q}\rfloor+1}}{1+q},\ \ \ \ 2-\sqrt{3}\le q<1.$$
Suppose $\frac{\sqrt{2q}}{1-q}\in [x,x+1)$ for some $x\in\mathbb{N}$. Then,
$f(q)=1-\frac{2q^{x+1}}{1+q}$ when $q\in[1+\frac{1-\sqrt{2x^2+1}}{x^2},1+\frac{1-\sqrt{2(x+1)^2+1}}{(x+1)^2})$. Since $f'(q)=-\frac{2q^x(x+1)+2q^{x+1}x}{(1+q)^2}<0$, $f(q)$ is decreasing in $[1+\frac{1-\sqrt{2x^2+1}}{x^2},1+\frac{1-\sqrt{2(x+1)^2+1}}{(x+1)^2})$. Thus,
$$\inf\limits_{q\in [1+\frac{1-\sqrt{2x^2+1}}{x^2},1+\frac{1-\sqrt{2(x+1)^2+1}}{(x+1)^2})}f(q)=\lim\limits_{q\uparrow 1+\frac{1-\sqrt{2(x+1)^2+1}}{(x+1)^2}}f(q)=1-\frac{2\left(1+\frac{1-\sqrt{2(x+1)^2+1}}{(x+1)^2}\right)^{x+1}}{2+\frac{1-\sqrt{2(x+1)^2+1}}{(x+1)^2}}.$$

Set $a_n:=1-\frac{2\left(1+\frac{1-\sqrt{2n^2+1}}{n^2}\right)^{n}}{2+\frac{1-\sqrt{2n^2+1}}{n^2}}=1-\frac{\left(1-\frac{2}{1+\sqrt{2n^2+1}}\right)^{n}}{1-\frac{1}{1+\sqrt{2n^2+1}}}$, $n\ge2$. Note that $\left(1-\frac{2}{1+\sqrt{2n^2+1}}\right)^{n}$ and $1-\frac{1}{1+\sqrt{2n^2+1}}$ are increasing with respect to $n$ for $n\ge2$. Then,
$$\frac{1}{4}\le\left(1-\frac{2}{1+\sqrt{2n^2+1}}\right)^{n}<\lim_{n\rightarrow\infty}\left(1-\frac{2}{1+\sqrt{2n^2+1}}\right)^{n}=e^{-\sqrt{2}},$$
and
$$\frac{3}{4}\le1-\frac{1}{1+\sqrt{2n^2+1}}<1.$$
Thus,
$$\frac{\left(1-\frac{2}{1+\sqrt{2n^2+1}}\right)^{n}}{1-\frac{1}{1+\sqrt{2n^2+1}}}<\frac{4}{3}e^{-\sqrt{2}}.$$
Hence,
$$a_n>1-\frac{4}{3}e^{-\sqrt{2}}\approx0.67584>\frac{\sqrt{3}}{3}.$$

By cases (i) and (ii), we get
$$
\inf_{p\in (0,1)}P\left\{|X^S_p-E[X^S_p]|\le \sqrt{{\rm Var}(X^S_p)}\right\}=\frac{\sqrt{3}}{3}.
$$
Therefore, the proof of (\ref{as2}) is complete by noting that
$$
P\left\{\left|X^S_{\sqrt{3}-1}-E\left[X^S_{\sqrt{3}-1}\right]\right|< \sqrt{{\rm Var}\left(X^S_{\sqrt{3}-1}\right)}\right\}=P\left\{\left|X^S_{\sqrt{3}-1}\right|< 1\right\}=\frac{\sqrt{3}}{3}.
$$
\hfill \fbox

\section{Poisson distribution and symmetric Poisson distribution}\setcounter{equation}{0}

\subsection{Poisson distribution}\label{41}
For $\lambda>0$, define
\begin{eqnarray*}
I_p(\lambda)=P\left\{|N_\lambda-E[N_\lambda]|\le \sqrt{{\rm Var}(N_\lambda)}\right\}.
\end{eqnarray*}

\begin{thm}\label{pro1c} We have
\begin{equation}\label{as3}
\inf_{\lambda>0}I_p(\lambda)=1.5e^{-1}\approx0.55182.
\end{equation}
\end{thm}

\noindent Proof.\ \  (i) For $\lambda\in(0,\frac{3-\sqrt{5}}{2}\approx0.38197)$, we have $\lambda-\sqrt{\lambda}<0$ and $0<\lambda+\sqrt{\lambda}<1$. Then,
\begin{eqnarray*}
I_p(\lambda)=P\{N_\lambda=0\}=e^{-\lambda}>e^{-\frac{3-\sqrt{5}}{2}}\approx0.68252.
\end{eqnarray*}

\noindent (ii)  For $\lambda\in[\frac{3-\sqrt{5}}{2},1)$, we have $\lambda-\sqrt{\lambda}<0$ and $1\le\lambda+\sqrt{\lambda}<2$. Then,
\begin{eqnarray*}
I_p(\lambda)=P\{N_\lambda=0\}+P\{N_\lambda=1\}=e^{-\lambda}(1+\lambda)>2e^{-1}\approx0.73576,
\end{eqnarray*}
where we have used the fact that $(e^{-\lambda}(1+\lambda))'=-\lambda e^{-\lambda}<0$.

\noindent (iii) For $\lambda=1$, we have
\begin{eqnarray*}
I_p(\lambda)=P\{N_\lambda=0\}+P\{N_\lambda=1\}+P\{N_\lambda=2\}=2.5e^{-1}\approx0.9197.
\end{eqnarray*}

\noindent (iv) For $\lambda\in(1,\frac{7-\sqrt{13}}{2}\approx1.69722)$, we have $0<\lambda-\sqrt{\lambda}<1$ and $2<\lambda+\sqrt{\lambda}<3$. Then,
\begin{eqnarray*}
I_p(\lambda)=P\{N_\lambda=1\}+P\{N_\lambda=2\}=e^{-\lambda}\left(\lambda+\frac{\lambda^2}{2}\right).
\end{eqnarray*}

Let $f_1(\lambda)=e^{-\lambda}(\lambda+\frac{\lambda^2}{2})$. Since  $f_1'(\lambda)=e^{-\lambda}(1-\frac{\lambda^2}{2})$,  $f_1(\lambda)$ is increasing in $(1,\sqrt{2})$  and decreasing in $(\sqrt{2},\frac{7-\sqrt{13}}{2})$. Note that $f_1(1)=1.5e^{-1}\approx0.55182$ and $f_1(\frac{7-\sqrt{13}}{2})\approx0.57476$. Then,
$$I_p(\lambda)>f_1(1)=1.5e^{-1}.$$

\noindent (v) For $\lambda\in[\frac{7-\sqrt{13}}{2},\frac{9-\sqrt{17}}{2})$, we have $0<\lambda-\sqrt{\lambda}<1$ and $3\le\lambda+\sqrt{\lambda}<4$. Then,
\begin{eqnarray*}
I_p(\lambda)=P\{N_\lambda=1\}+P\{N_\lambda=2\}+P\{N_\lambda=3\}=e^{-\lambda}\left(\lambda+\frac{\lambda^2}{2}+\frac{\lambda^3}{6}\right).
\end{eqnarray*}

Let $f_2(\lambda)=e^{-\lambda}(\lambda+\frac{\lambda^2}{2}+\frac{\lambda^3}{6})$. Since $f'_2(\lambda)=e^{-\lambda}(1-\frac{\lambda^3}{6})$,  $f_2(\lambda)$  is increasing in $[\frac{7-\sqrt{13}}{2},\sqrt[3]{6}\approx1.8171]$  and decreasing in $(\sqrt[3]{6},\frac{9-\sqrt{17}}{2})$. Note that $f_2(\frac{9-\sqrt{17}}{2})\in(0.68335,0.68336)$ and $ f_2(\frac{7-\sqrt{13}}{2})\approx 0.72403$. Then,
$$I_p(\lambda)>f_2\left(\frac{9-\sqrt{17}}{2}\right)>0.68335.$$

\noindent (vi) For $\lambda\in[\frac{9-\sqrt{17}}{2},\frac{3+\sqrt{5}}{2}\approx2.618]$, we have $0<\lambda-\sqrt{\lambda}\le1$ and $4\le\lambda+\sqrt{\lambda}<5$. Then,
\begin{eqnarray*}
I_p(\lambda)=P\{N_\lambda=1\}+P\{N_\lambda=2\}+P\{N_\lambda=3\}+P\{N_\lambda=4\}=e^{-\lambda}\left(\lambda+\frac{\lambda^2}{2}+\frac{\lambda^3}{6}+\frac{\lambda^4}{24}\right).
\end{eqnarray*}

Let $f_3(\lambda)=e^{-\lambda}(\lambda+\frac{\lambda^2}{2}+\frac{\lambda^3}{6}+\frac{\lambda^4}{24})$. Since $f'_3(\lambda)=e^{-\lambda}(1-\frac{\lambda^4}{24})$,  $f_3(\lambda)$ is decreasing in $[\frac{9-\sqrt{17}}{2},\frac{3+\sqrt{5}}{2}\approx2.618]$. Note that $f_3(\frac{3+\sqrt{5}}{2})\in(0.80191,0.80192)$. Then,
$$I_p(\lambda)>0.80191.$$

\noindent (vii)  For $\lambda\in(\frac{3+\sqrt{5}}{2},\frac{11-\sqrt{21}}{2}\approx3.20871)$, we have $1<\lambda-\sqrt{\lambda}<2$ and $4<\lambda+\sqrt{\lambda}<5$. Then,
\begin{eqnarray*}
I_p(\lambda)=P\{N_\lambda=2\}+P\{N_\lambda=3\}+P\{N_\lambda=4\}=e^{-\lambda}\left(\frac{\lambda^2}{2}+\frac{\lambda^3}{6}+\frac{\lambda^4}{24}\right).
\end{eqnarray*}

Let $f_4(\lambda)=e^{-\lambda}(\frac{\lambda^2}{2}+\frac{\lambda^3}{6}+\frac{\lambda^4}{24})$. Since $f'_4(\lambda)=e^{-\lambda}(\lambda-\frac{\lambda^4}{24})$, $f_4(\lambda)$ is increasing in $(\frac{3+\sqrt{5}}{2},\sqrt[3]{24}\approx2.8845]$  and decreasing in $(\sqrt[3]{24},\frac{11-\sqrt{21}}{2})$. Note that $f_4(\frac{3+\sqrt{5}}{2})\approx0.61094$ and $f_4(\frac{11-\sqrt{21}}{2})\in(0.60899,0.609)$. Then,
$$I_p(\lambda)>0.60899.$$

\noindent (viii) For $\lambda\in[\frac{11-\sqrt{21}}{2},4)$, we have $1<\lambda-\sqrt{\lambda}<2$ and $5\le\lambda+\sqrt{\lambda}<6$. Then,
\begin{eqnarray*}
I_p(\lambda)&=&P\{N_\lambda=2\}+P\{N_\lambda=3\}+P\{N_\lambda=4\}+P\{N_\lambda=5\}\\
&=&e^{-\lambda}\left(\frac{\lambda^2}{2}+\frac{\lambda^3}{6}+\frac{\lambda^4}{24}+\frac{\lambda^5}{120}\right).
\end{eqnarray*}

Let $f_5(\lambda)=e^{-\lambda}(\frac{\lambda^2}{2}+\frac{\lambda^3}{6}+\frac{\lambda^4}{24}+\frac{\lambda^5}{120})$. Since $f'_5(\lambda)=e^{-\lambda}(\lambda-\frac{\lambda^5}{120})$,  $f_5(\lambda)$ is increasing in $[\frac{11-\sqrt{21}}{2},\sqrt[4]{120}\approx3.30975]$  and decreasing in $(\sqrt[4]{120},4)$. Note that $f_5(4)\in(0.69355,0.69356)$ and $f_5(\frac{11-\sqrt{21}}{2})\approx0.72353$. Then,
$$I_p(\lambda)>0.69355.$$

\noindent (ix) For $\lambda=4$, we have $\lambda-\sqrt{\lambda}=2$ and $\lambda+\sqrt{\lambda}=6$. Then,
\begin{eqnarray*}
I_p(\lambda)&=&P\{N_\lambda=2\}+P\{N_\lambda=3\}+P\{N_\lambda=4\}+P\{N_\lambda=5\}+P\{N_\lambda=6\}\\
&=&e^{-\lambda}\left(\frac{\lambda^2}{2}+\frac{\lambda^3}{6}+\frac{\lambda^4}{24}+\frac{\lambda^5}{120}+\frac{\lambda^6}{720}\right).
\end{eqnarray*}
Thus,
$$I_p(4)>0.79774.$$

\noindent (x) For $\lambda\in(4, \frac{15-\sqrt{29}}{2}\approx4.80742)$, we have $2<\lambda-\sqrt{\lambda}<3$ and $6<\lambda+\sqrt{\lambda}<7$. Then,
\begin{eqnarray*}
I_p(\lambda)&=&P\{N_\lambda=3\}+P\{N_\lambda=4\}+P\{N_\lambda=5\}+P\{N_\lambda=6\}\\
&=&e^{-\lambda}\left(\frac{\lambda^3}{6}+\frac{\lambda^4}{24}+\frac{\lambda^5}{120}+\frac{\lambda^6}{720}\right).
\end{eqnarray*}

Let $f_6(\lambda)=e^{-\lambda}(\frac{\lambda^3}{6}+\frac{\lambda^4}{24}+\frac{\lambda^5}{120}+\frac{\lambda^6}{720})$. Since $f'_6(\lambda)=e^{-\lambda}(\frac{\lambda^2}{2}-\frac{\lambda^6}{720})$,  $f_6(\lambda)$ is decreasing in $(4, \frac{15-\sqrt{29}}{2})$ and $f_6(\frac{15-\sqrt{29}}{2})\in(0.64792,0.64793)$. Then,
$$I_p(\lambda)>0.64792.$$

\noindent (xi) For $\lambda\ge \frac{15-\sqrt{29}}{2}$, let $g_1(\lambda):=P\{N_\lambda\le \lambda+\sqrt{\lambda}\}$ and $g_2(\lambda):=P\{N_\lambda< \lambda-\sqrt{\lambda}\}$. Then,
$$I_p(\lambda)=g_1(\lambda)-g_2(\lambda).$$

\noindent (xia) Suppose $\lambda=n$ for some $n\in\mathbb{N}$. Let $Y_1,Y_2,\dots,Y_n$ be i.i.d. Poisson(1) random variables. We have $E[Y_1-1]=0$, ${\rm Var}(Y_1-1)=1$ and
$$N_n=Y_1+Y_2+\cdots+Y_n\ \ {\rm in\ distribution}.$$
Denote by $\Phi$ the cumulative distribution function of the standard normal random variable. By the Berry-Esseen theorem, we get
 $$\left|P\left\{\frac{1}{\sqrt{n}}(Y_1-1+\cdots+Y_n-1)\le1\right\}-\Phi(1)\right| \le \frac{C\rho}{\sqrt{n}}, $$
 where $\rho=E[|Y_1-1|^3]=E[(Y_1-1)^3]+2P\{Y_1=0\}=E[Y_1^3]-3E[Y_1^2]+3E[Y_1]-1+2e^{-1}=1+2e^{-1}$ and   $C$  can be taken to 0.7655 (cf. \cite{Sh}). We have
 \begin{eqnarray}\label{K1}
  & &\left|P\left\{\frac{1}{\sqrt{n}}(Y_1-1+\cdots+Y_n-1)\le1\right\}-\Phi(1)\right| \le \frac{C\rho}{\sqrt{n}}\nonumber\\
    &\Leftrightarrow&\left|P\left\{N_n\le n+\sqrt{n}\right\}-\Phi(1)\right| \le \frac{C\rho}{\sqrt{n}}\nonumber\\
     &\Leftrightarrow&\Phi(1)-\frac{C\rho}{\sqrt{n}}\le P\left\{N_n\le n+\sqrt{n}\right\} \le \Phi(1)+\frac{C\rho}{\sqrt{n}}.
 \end{eqnarray}

Note that $C\rho<0.7656(1+2e^{-1})<1.328898$. Then, for any $ n\ge604$, we have
 $$P\left\{N_n\le n+\sqrt{n}\right\}\ge \Phi(1)-\frac{1.328898}{\sqrt{604}}>0.7872.$$

\noindent (xib) Consider $\lambda\in (n,n+1)$ for $n\ge604$. We have
 $$P\{N_\lambda\le \lambda+\sqrt{\lambda}\}\ge P\{N_\lambda\le n+\sqrt{n}\}>P\{N_{n+1}\le n+\sqrt{n}\},$$
and
  \begin{eqnarray}\label{K2}
  & &\left|P\left\{\frac{1}{\sqrt{n+1}}(Y_1-1+\cdots+Y_n-1+Y_{n+1}-1)\le\frac{\sqrt{n}-1}{\sqrt{n+1}}\right\}-\Phi\left(\frac{\sqrt{n}-1}{\sqrt{n+1}}\right)\right| \le \frac{C\rho}{\sqrt{n+1}} \nonumber \\
    &\Leftrightarrow&\left|P\left\{\frac{N_{n+1}-n-1}{\sqrt{n+1}}\le\frac{\sqrt{n}-1}{\sqrt{n+1}}\right\}-\Phi\left(\frac{\sqrt{n}-1}{\sqrt{n+1}}\right)\right| \le\frac{C\rho}{\sqrt{n+1}}\nonumber\\
    &\Leftrightarrow&\left|P\left\{N_{n+1}\le n+\sqrt{n}\right\}-\Phi\left(\frac{\sqrt{n}-1}{\sqrt{n+1}}\right)\right| \le \frac{C\rho}{\sqrt{n+1}}\nonumber\\
     &\Leftrightarrow&\Phi\left(\frac{\sqrt{n}-1}{\sqrt{n+1}}\right)-\frac{C\rho}{\sqrt{n+1}}\le P\left\{N_{n+1}\le n+\sqrt{n}\right\} \le \Phi\left(\frac{\sqrt{n}-1}{\sqrt{n+1}}\right)+\frac{C\rho}{\sqrt{n+1}}.
 \end{eqnarray}
Since $\Phi(\frac{\sqrt{n}-1}{\sqrt{n+1}})-\frac{C\rho}{\sqrt{n+1}}=\Phi(\sqrt{1-\frac{1}{n+1}}-\sqrt{\frac{1}{n+1}})-\frac{C\rho}{\sqrt{n+1}}$ is increasing with respect to $n$, we get
 \begin{eqnarray*}
P\left\{N_{n+1}\le n+\sqrt{n}\right\}&\ge& \Phi\left(\frac{\sqrt{n}-1}{\sqrt{n+1}}\right)-\frac{C\rho}{\sqrt{n+1}}\\
&\ge & \Phi\left(\frac{\sqrt{604}-1}{\sqrt{604+1}}\right)-\frac{C\rho}{\sqrt{604+1}}\\
 &>&\Phi(0.9585)-0.0541\\
&>&0.7769.
 \end{eqnarray*}
Then,
\begin{equation}\label{g11}
 g_1(\lambda)>0.7769,\ \ \ \ \forall \lambda\ge604.
\end{equation}

\noindent (xic) We have
$$g_2(\lambda)\le P\{N_\lambda\le \lambda-\sqrt{\lambda}\}. $$
 Suppose $\lambda=n$ for some $n\in\mathbb{N}$. By the Berry-Esseen theorem, we get
 \begin{eqnarray}\label{K3}
  & &\left|P\left\{\frac{1}{\sqrt{n}}(Y_1-1+\cdots+Y_n-1)\le-1\right\}-\Phi(-1)\right| \le \frac{C\rho}{\sqrt{n}} \nonumber \\
    &\Leftrightarrow&\left|P\left\{\frac{N_n-n}{\sqrt{n}}\le-1\right\}-\Phi(-1)\right| \le \frac{C\rho}{\sqrt{n}} \nonumber\\
    &\Leftrightarrow&\left|P\left\{N_n\le n-\sqrt{n}\right\}-\Phi(-1)\right| \le \frac{C\rho}{\sqrt{n}}\nonumber\\
     &\Leftrightarrow&\Phi(-1)-\frac{C\rho}{\sqrt{n}}\le P\left\{N_n\le n-\sqrt{n}\right\} \le \Phi(-1)+\frac{C\rho}{\sqrt{n}}.
 \end{eqnarray}
Then, for any $n\ge604$, we have
 $$P\left\{N_n\le n-\sqrt{n}\right\}\le  \Phi(-1)+\frac{1.328898}{\sqrt{604}}<0.2128.$$

\noindent (xid) Consider $\lambda\in (n,n+1)$ for $n\ge604$.  We have
 $$P\{N_\lambda\le \lambda-\sqrt{\lambda}\}< P\{N_n\le \lambda-\sqrt{\lambda}\}\le P\{N_{n}\le n+1-\sqrt{n+1}\}.$$
By the Berry-Esseen theorem, we get
  \begin{eqnarray}\label{K4}
  & &\left|P\left\{\frac{1}{\sqrt{n}}(Y_1-1+\cdots+Y_{n}-1)\le\frac{1-\sqrt{n+1}}{\sqrt{n}}\right\}-\Phi\left(\frac{1-\sqrt{n+1}}{\sqrt{n}}\right)\right| \le \frac{C\rho}{\sqrt{n}} \nonumber \\
    &\Leftrightarrow&\left|P\left\{\frac{N_{n}-n}{\sqrt{n}}\le\frac{1-\sqrt{n+1}}{\sqrt{n}}\right\}-\Phi\left(\frac{1-\sqrt{n+1}}{\sqrt{n}}\right)\right| \le \frac{C\rho}{\sqrt{n}}\nonumber\\
    &\Leftrightarrow&\left|P\left\{N_{n}\le n+1-\sqrt{n+1}\right\}-\Phi\left(\frac{1-\sqrt{n+1}}{\sqrt{n}}\right)\right| \le \frac{C\rho}{\sqrt{n}}\nonumber\\
     &\Leftrightarrow&\Phi\left(\frac{1-\sqrt{n+1}}{\sqrt{n}}\right)-\frac{C\rho}{\sqrt{n}}\le P\left\{N_{n}\le n+1-\sqrt{n+1}\right\} \le \Phi\left(\frac{1-\sqrt{n+1}}{\sqrt{n}}\right)+\frac{C\rho}{\sqrt{n}}.\ \ \ \
 \end{eqnarray}

Since $ \Phi(\frac{1-\sqrt{n+1}}{\sqrt{n}})+\frac{C\rho}{\sqrt{n}}$ is decreasing with respect to $n$, we get
\begin{eqnarray*}
 P\left\{N_{n}\le n+1-\sqrt{n+1}\right\} &\le&  \Phi\left(\frac{1-\sqrt{604+1}}{\sqrt{604}}\right)+\frac{C\rho}{\sqrt{604}}\\
&<&\Phi(-0.96)+0.0541\\
&<&0.2227.
\end{eqnarray*}
 Then,
\begin{equation}\label{g22}
  g_2(\lambda)<0.2227, \ \ \ \ \forall \lambda\ge604.
\end{equation}
Thus, by (\ref{g11}) and (\ref{g22}),  we find that for $\lambda\ge604$,
 $$I_p(\lambda)=g_1(\lambda)-g_2(\lambda)>0.7769-0.2227=0.5542>1.5e^{-1}.$$

\noindent (xie) Finally, we consider the case that  $\lambda\in[\frac{15-\sqrt{29}}{2},604)$.

\noindent (xie1) For $n\in\mathbb{N}$, let $\lambda_{n,1}^{(1)},\lambda_{n,2}^{(1)}$ be the solutions of $\lambda+\sqrt{\lambda}=n$, where
$$\lambda_{n,1}^{(1)}=\frac{2n+1-\sqrt{4n+1}}{2},\ \ \ \ \lambda_{n,2}^{(1)}=\frac{2n+1+\sqrt{4n+1}}{2}.$$
Since $\lambda<n$, the equation has a unique solution, which is denoted by $ \lambda_{n}^{(1)}=\frac{2n+1-\sqrt{4n+1}}{2}<n-1$.

For fixed $x$, $P\{N_\lambda\le x\}$ is decreasing with respect to $\lambda$. Then,  when $\lambda+\sqrt{\lambda}\in[n,n+1)$, $P\{N_\lambda\le\lambda+\sqrt{\lambda}\}=P\{N_\lambda\le n\}$  is decreasing for $\lambda\in[\lambda_{n}^{(1)},\lambda_{n+1}^{(1)})$. Thus,  for $\lambda\in[\lambda_{n}^{(1)},\lambda_{n+1}^{(1)})$, we have
$$g_1(\lambda)> P\{N_{\lambda_{n+1}^{(1)}}\le n\}.$$
Note that $\frac{15-\sqrt{29}}{2}=\lambda_{7}^{(1)}$ and $ 604\in(\lambda_{628}^{(1)},\lambda_{629}^{(1)})$. For $\lambda\in[\frac{15-\sqrt{29}}{2},604)$,  we use Matlab to find the minimal value of $P\{N_{\lambda_{n}^{(1)}}\le n-1\}$ for $ n\in[8,629]$, which is $P\{N_{\lambda_{8}^{(1)}}\le 7\}\in(0.79345,0.79346)$. See Figure 1. We refer to the Appendix for Matlab codes. Hence,  we have
\begin{equation}\label{g1}
  g_1(\lambda)>0.79345,\ \ \ \ \lambda\in\left[\frac{15-\sqrt{29}}{2},604\right).
\end{equation}
\begin{figure}[!h]
   \centering
   \includegraphics[width=8cm,height=5cm]{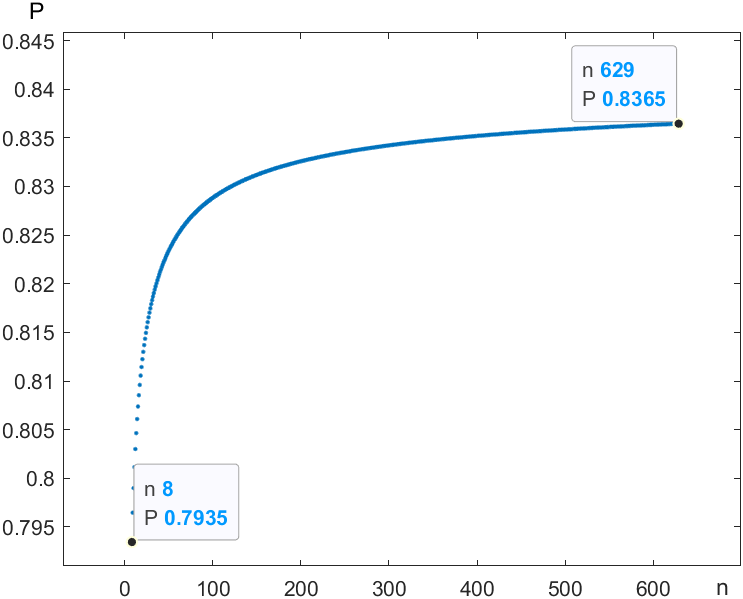}
    \caption{Figure of $P\{X_{\lambda_{n}^{(1)}}\le n-1\}$ for $n\in[8,629]$.}
\end{figure}

\noindent (xie2) For $n\in\mathbb{N}$, let $\lambda_{n,1}^{(2)},\lambda_{n,2}^{(2)}$ be the solutions of $\lambda-\sqrt{\lambda}=n$, where
$$\lambda_{n,1}^{(2)}=\frac{2n+1-\sqrt{4n+1}}{2},\ \ \ \ \lambda_{n,2}^{(2)}=\frac{2n+1+\sqrt{4n+1}}{2}.$$
Since $\lambda>n$, the equation has a unique solution, which is denoted by $\lambda_{n}^{(2)}=\frac{2n+1+\sqrt{4n+1}}{2}>n+1$.

For fixed $x$, $P\{N_\lambda< x\}$ is decreasing with respect to $\lambda$. Then, when $\lambda-\sqrt{\lambda}\in(n,n+1]$, $P\{N_\lambda<\lambda-\sqrt{\lambda}\}=P\{N_\lambda\le n\}$ is decreasing for $\lambda\in(\lambda_{n}^{(2)},\lambda_{n+1}^{(2)}]$. Thus, for $\lambda\in (\lambda_{n}^{(2)},\lambda_{n+1}^{(2)}]$, we have
$$g_2(\lambda)< P\{N_{\lambda_{n}^{(2)}}\le n\}.$$
Note that $\frac{15-\sqrt{29}}{2}\in(\lambda_{2}^{(2)},\lambda_{3}^{(2)})$ and $604\in(\lambda_{579}^{(2)},\lambda_{580}^{(2)})$. For $\lambda\in[\frac{15-\sqrt{29}}{2},604)$, we consider the maximal value of $g_2(\frac{15-\sqrt{29}}{2})$ and  $P\{N_{\lambda_{n}^{(2)}}\le n\}$ for $ n\in[3,579]$. We have $g_2(\frac{15-\sqrt{29}}{2})=P\{N_{\frac{15-\sqrt{29}}{2}}\le2\}\approx0.14184$. By virtue of Matlab, we find the maximal value of  $P\{N_{\lambda_{n}^{(2)}}\le n\}$ for $ n\in[3,579]$, which is $P\{N_{\lambda_{3}^{(2)}}\le 3\}\in(0.22506,0.22507)$. See Figure 2. We refer to the Appendix for  Matlab codes.  Hence, we have
\begin{equation}\label{g2}
  g_2(\lambda)<0.22507,\ \ \ \ \lambda\in\left[\frac{15-\sqrt{29}}{2},604\right).
\end{equation}
\newpage
\begin{figure}[!h]
    \centering
   \includegraphics[width=8cm,height=5cm]{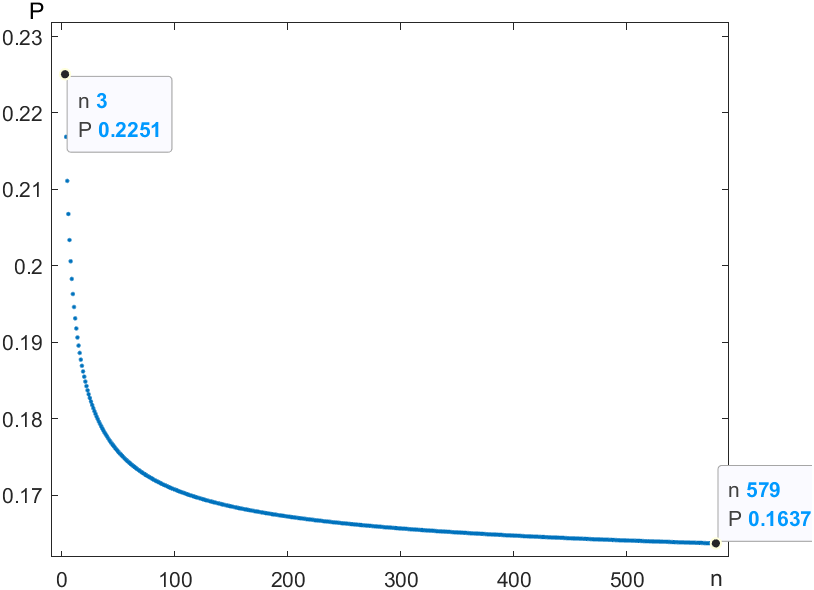}
    \caption{Figure of $P\{X_{\lambda_{n}^{(2)}}\le n\}$ for $n\in[3,579]$.}
\end{figure}

By (\ref{g1}) and (\ref{g2}), we conclude that for $\lambda\in[\frac{15-\sqrt{29}}{2},604)$,
$$
  I_p(\lambda)=g_1(\lambda)-g_2(\lambda)>0.79345-0.22507=0.56838>1.5e^{-1}.
$$
Therefore, (\ref{as3}) holds based on all the above cases. \hfill \fbox

\begin{thm}\label{pro2} We have
$$
\inf_{\lambda>0}P\left\{|N_\lambda-E[N_\lambda]|< \sqrt{{\rm Var}(N_\lambda)}\right\}=P\left\{|N_1-E[N_1]|< \sqrt{{\rm Var}(N_1)}\right\}=e^{-1},
$$
which implies that
$$
P_{{\cal I}_0}\le e^{-1}\approx 0.36788.
$$
\end{thm}

Using an argument similar to the proof of Theorem \ref{pro1c}, we can prove Theorem \ref{pro2}. The details will be given in the Appendix.

\subsection{Symmetric Poisson distribution}\label{42}

Let $N_{\lambda}$ and $N'_{\lambda}$ be  independent Poisson random variables with the same parameter $\lambda>0$. Define
$$
Y_{\lambda}=N_{\lambda}-N'_{\lambda},
$$
and
$$
P_{\lambda}:=P\left\{|Y_{\lambda}-E[Y_{\lambda}]|\le \sqrt{{\rm Var}(Y_{\lambda})}\right\},\ \ \ \ P^0_{\lambda}:=P\left\{|Y_{\lambda}-E[Y_{\lambda}]|< \sqrt{{\rm Var}(Y_{\lambda})}\right\}.
$$
Then,
$$E[Y_\lambda]=0,\ \ \ \ {\rm Var}(Y_\lambda)=2\lambda.$$
\begin{thm}\label{thm2} We have
\begin{eqnarray*}
\inf\limits_{\lambda>0} P_{\lambda}=\inf\limits_{\lambda>0} P^0_{\lambda}=P^0_{\frac{1}{2}}=e^{-1}\sum_{k=0}^{\infty}\frac{1}{2^{2k}(k!)^2},
\end{eqnarray*}
which implies that
$$
P_{\cal I}\le e^{-1}\sum_{k=0}^{\infty}\frac{1}{2^{2k}(k!)^2}\approx 0.46576.
$$
\end{thm}

\noindent Proof.\ \ For $\lambda\in(\frac{(n-1)^2}{2},\frac{n^2}{2})$, $n\in\mathbb{N}$, we have
\begin{eqnarray}\label{10a}
P_{\lambda}=P^0_{\lambda}&=&P\left\{N_{\lambda}-N'_{\lambda}=0\right\}+2\sum_{l=1}^{n-1}P\left\{N_{\lambda}-N'_{\lambda}=l\right\}\nonumber\\
&=&e^{-2\lambda}\left[\sum_{k=0}^{\infty}\frac{\lambda^{2k}}{(k!)^2}+2\sum_{l=1}^{n-1}\sum_{k=0}^{\infty}\frac{\lambda^{2k+l}}{k!(k+l)!}\right].
\end{eqnarray}
For $\lambda=\frac{(n-1)^2}{2}$, $n\ge2$, we have
\begin{eqnarray}\label{10a-1}
P_{\lambda}&=&P\left\{N_{\lambda}-N'_{\lambda}=0\right\}+2\sum_{l=1}^{n-1}P\left\{N_{\lambda}-N'_{\lambda}=l\right\}\nonumber\\
&=&e^{-2\lambda}\left[\sum_{k=0}^{\infty}\frac{\lambda^{2k}}{(k!)^2}+2\sum_{l=1}^{n-1}\sum_{k=0}^{\infty}\frac{\lambda^{2k+l}}{k!(k+l)!}\right].
\end{eqnarray}
For $\lambda=\frac{n^2}{2}$, $n\in\mathbb{N}$, we have
\begin{eqnarray}\label{10a-2}
P^0_{\lambda}&=&P\left\{N_{\lambda}-N'_{\lambda}=0\right\}+2\sum_{l=1}^{n-1}P\left\{N_{\lambda}-N'_{\lambda}=l\right\}\nonumber\\
&=&e^{-2\lambda}\left[\sum_{k=0}^{\infty}\frac{\lambda^{2k}}{(k!)^2}+2\sum_{l=1}^{n-1}\sum_{k=0}^{\infty}\frac{\lambda^{2k+l}}{k!(k+l)!}\right].
\end{eqnarray}
By (\ref{10a})-(\ref{10a-2}), we know that the function $\lambda\mapsto P_{\lambda}$ is right continuous in $(0,\infty)$ and the function $\lambda\mapsto P^0_{\lambda}$ is left continuous in $(0,\infty)$. Below are graphs of  $P_{\lambda}$ and $P^0_{\lambda}$ for $\lambda\in(0,\frac{1}{2})$ and $\lambda\in(0,10)$, respectively.
\begin{figure}[h]
\begin{center}
\scalebox{0.5}{\includegraphics{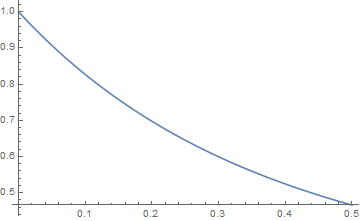}}
\end{center}
\end{figure}
\begin{center}
{\small Figure 3: Functions $P_{\lambda}$  and $P^0_{\lambda}$ for $\lambda\in(0,\frac{1}{2})$.}
\end{center}

\newpage
\begin{figure}[h]
\begin{center}
\scalebox{0.5}{\includegraphics{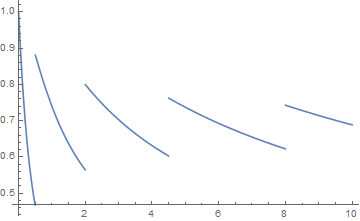}}
\end{center}
\end{figure}
\begin{center}
{\small Figure 4: Functions $P_{\lambda}$  and $P^0_{\lambda}$ for $\lambda\in(0,10)$.}
\end{center}

\noindent {\it Step 1:} Consider $\lambda\in(0,\frac{1}{2})$. We have
\begin{eqnarray*}
P^0_{\lambda}=P\left\{N_{\lambda}-N'_{\lambda}=0\right\}=e^{-2\lambda}\sum_{k=0}^{\infty}\frac{\lambda^{2k}}{(k!)^2}.
\end{eqnarray*}
Then,
\begin{eqnarray*}
\frac{dP^0_{\lambda}}{d\lambda}&=&e^{-2\lambda}\left[(-2)+(-2)\sum_{k=1}^{\infty}\frac{\lambda^{2k}}{(k!)^2}+\sum_{k=1}^{\infty}\frac{2k\lambda^{2k-1}}{(k!)^2}\right]\\
&<&2e^{-2\lambda}\left[-1+\sum_{k=1}^{\infty}\frac{k\lambda^{2k-1}}{(k!)^2}\right]\\
&<&2e^{-2\lambda}\left[-1+\lambda+\frac{\lambda^3}{2}+\sum_{k=3}^{\infty}\frac{\lambda^{2k-1}}{k!}\right]\\
&<&2e^{-2\lambda}\left[-1+\lambda+\frac{\lambda^3}{2}+\lambda^2e^{\lambda}\right]\\
&<&2e^{-2\lambda}\left[-1+0.5+\frac{(0.5)^3}{2}+(0.5)^2e^{0.5}\right]\\
&\approx&2e^{-2\lambda}(-0.0253)\\
&<&0.
\end{eqnarray*}
Thus, the function $\lambda\mapsto P^0_{\lambda}$ is decreasing in $(0,\frac{1}{2})$. Hence,
$$
\inf\limits_{\lambda\in(0,\frac{1}{2})} P_{\lambda}=\inf\limits_{\lambda\in(0,\frac{1}{2})} P^0_{\lambda}=P^0_{\frac{1}{2}}=e^{-1}\sum_{k=0}^{\infty}\frac{1}{2^{2k}(k!)^2}.
$$

\noindent {\it Step 2:} Consider $\lambda\in[200,\infty)$. Set $Y_0=0$. We have
$$
Y_{\lambda}=\sum_{k=1}^{\lfloor \lambda\rfloor}\left(Y_{\frac{k\lambda}{\lfloor \lambda\rfloor}}-Y_{\frac{(k-1)\lambda}{\lfloor \lambda\rfloor}}\right).
$$
Define
$$
\rho=E\left[\left|Y_{\frac{\lambda}{\lfloor \lambda\rfloor}}\right|^3\right].
$$
By the Berry-Esseen theorem, we get
\begin{eqnarray}\label{91}
P^0_{\lambda}
&=&P\left\{|Y_{\lambda}|< \sqrt{2\lambda}\right\}\nonumber\\
&=&P\left\{Y_{\lambda}< \sqrt{2\lambda}\right\}-P\left\{Y_{\lambda}\le -\sqrt{2\lambda}\right\}\nonumber\\
&=&P\left\{\frac{Y_{\lambda}}{\sqrt{2\lambda}}< 1\right\}-P\left\{\frac{Y_{\lambda}}{\sqrt{2\lambda}}\le -1\right\}\nonumber\\
&\ge&\Phi(1)-\Phi(-1)-\frac{2C\rho}{\left(\frac{2\lambda}{\lfloor \lambda\rfloor}\right)^{\frac{3}{2}}\sqrt{\lfloor \lambda\rfloor}}.
\end{eqnarray}

Note that
$$
E[N^3_{\lambda}]=\lambda(1+3\lambda+\lambda^2),\ \ \ \ E[N^4_{\lambda}]=\lambda(1+7\lambda+6\lambda^2+\lambda^3).
$$
Then, for $\lambda\ge 200$, we have
\begin{eqnarray}\label{92}
\rho&=&E\left[\left|Y_{\frac{\lambda}{\lfloor \lambda\rfloor}}\right|^3\right]\nonumber\\
&\le&\left\{E\left[Y_{\frac{\lambda}{\lfloor \lambda\rfloor}}^2\right]\right\}^{\frac{1}{2}}\left\{E\left[Y_{\frac{\lambda}{\lfloor \lambda\rfloor}}^4\right]\right\}^{\frac{1}{2}}\nonumber\\
&=&\left\{{\rm Var}\left(Y_{\frac{\lambda}{\lfloor \lambda\rfloor}}\right)\right\}^{\frac{1}{2}}\left\{2E\left[N_{\frac{\lambda}{\lfloor \lambda\rfloor}}^4\right]-8E\left[N_{\frac{\lambda}{\lfloor \lambda\rfloor}}\right]E\left[N_{\frac{\lambda}{\lfloor \lambda\rfloor}}^3\right]+6\left(E\left[N^2_{\frac{\lambda}{\lfloor \lambda\rfloor}}\right]\right)^2\right\}^{\frac{1}{2}}\nonumber\\
&=&\left\{\frac{2\lambda}{\lfloor \lambda\rfloor}\right\}^{\frac{1}{2}}\left\{\frac{2\lambda}{\lfloor \lambda\rfloor}\left(1+\frac{6\lambda}{\lfloor \lambda\rfloor}\right)\right\}^{\frac{1}{2}}\nonumber\\
&<&\frac{2\lambda}{\lfloor \lambda\rfloor}\left(7+\frac{6}{200}\right)^{\frac{1}{2}}\nonumber\\
&<&\frac{5.31\lambda}{\lfloor \lambda\rfloor}.
\end{eqnarray}
Thus, by (\ref{91}) and (\ref{92}), we get
\begin{eqnarray*}
P^0_{\lambda}
>0.6826-\frac{2\cdot 0.7656\cdot 5.31}{2^{\frac{3}{2}} \sqrt{200}}>0.4793.
\end{eqnarray*}

\noindent {\it Step 3:} Consider $\lambda\in(\frac{1}{2},200)$. For $\lambda\in\left(\frac{(n-1)^2}{2},\frac{n^2}{2}\right)$ with $2\le n\le 20$, by (\ref{10a}), we get
\begin{eqnarray}\label{94}
\frac{dP^0_{\lambda}}{d\lambda}&=&-2e^{-2\lambda}\left[\sum_{k=0}^{\infty}\frac{\lambda^{2k}}{(k!)^2}+2\sum_{l=1}^{n-1}\sum_{k=0}^{\infty}\frac{\lambda^{2k+l}}{k!(k+l)!}\right]\nonumber\\
&&+e^{-2\lambda}\left[\sum_{k=1}^{\infty}\frac{2k\lambda^{2k-1}}{(k!)^2}+2\sum_{l=1}^{n-1}\sum_{k=0}^{\infty}\frac{(2k+l)\lambda^{2k+l-1}}{k!(k+l)!}\right]\nonumber\\
&>&-2e^{-2\lambda}\left[\sum_{k=0}^{\infty}\frac{\lambda^{2k}}{(k!)^2}+2\sum_{l=1}^{n-1}\sum_{k=0}^{\infty}\frac{\lambda^{2k+l}}{k!(k+l)!}\right]\nonumber\\
&>&-2e^{-2\lambda}\left[\left\{\sum_{k=0}^{\infty}\frac{\lambda^{k}}{k!}\right\}^2+2(n-1)\left\{\sum_{k=0}^{\infty}\frac{\lambda^{k}}{k!}\right\}^2\right]\nonumber\\
&=&-2(2n-1)\nonumber\\
&\ge&-78.
\end{eqnarray}

For $\lambda\in\left[\frac{1}{2},200\right]$, define
\begin{eqnarray}\label{97}
Q_{\lambda}=e^{-2\lambda}\left[\sum_{k=0}^{250}\frac{\lambda^{2k}}{(k!)^2}+2\sum_{l=1}^{\lfloor \sqrt{2\lambda}\rfloor}\sum_{k=0}^{250}\frac{\lambda^{2k+l}}{k!(k+l)!}\right].
\end{eqnarray}
By virtue of Mathematica, we find that
\begin{equation}\label{95}
\min\left\{Q_{\frac{1}{2}+\frac{r}{2\cdot 10^{3}}}:0\le r\le 399\cdot10^3\right\}=0.564565.
\end{equation}
We refer to the Appendix for Mathematica codes.

Suppoe $\lambda\in\left(\frac{(n-1)^2}{2}+\frac{r-1}{2\cdot 10^{3}},\frac{(n-1)^2}{2}+\frac{r}{2\cdot 10^{3}}\right)$ for some $1\le r\le (2n-1)10^3$ and $2\le n\le 20$.  Then, by the mean value theorem, (\ref{10a}) and  (\ref{94})--(\ref{95}), we get
\begin{eqnarray*}
P^0_{\lambda}&>&\left.e^{-2\lambda}\left[\sum_{k=0}^{\infty}\frac{\lambda^{2k}}{(k!)^2}+2\sum_{l=1}^{n-1}\sum_{k=0}^{\infty}\frac{\lambda^{2k+l}}{k!(k+l)!}\right]\right|_{\lambda=\frac{(n-1)^2}{2}+\frac{r-1}{2\cdot 10^{3}}}-78\cdot \frac{1}{2\cdot 10^{3}}\\
&>& Q_{\frac{(n-1)^2}{2}+\frac{r-1}{2\cdot 10^{3}}}-0.039\\
&>&0.56456-0.039\\
&=&0.52556.
\end{eqnarray*}
Therefore, the proof is complete based on all the above cases. \hfill \fbox

\section{Appendix}\setcounter{equation}{0}

\subsection{Matlab codes for \S \ref{41}}

(1) The minimal value of $P\{N_{\lambda_n^{(1)}}\le n-1\}$ for $[8,629]$:

clc

clear

n=8:1:629;

for i=1:1:length(n)

lam(i)=(2*n(i)+1-sqrt(4*n(i)+1))/2;

y(i)=poisscdf(n(i)-1,lam(i));

end

plot(n,y,'.')

min=min(y)
\vskip 0.3cm
Output result: min= 0.793450747058153
\vskip 0.4cm
(2) The maximal value of $P\{N_{\lambda_n^{(2)}}\le n\}$ for [3, 579]:

clc

clear

n=3:1:579;

for i=1:1:length(n)

lam(i)=(2*n(i)+1+sqrt(4*n(i)+1))/2;

y(i)=poisscdf(n(i),lam(i));

end

plot(n,y,'.')

max=max(y)
\vskip 0.3cm
Output result: max= 0.225065994481669

\subsection{Mathematica codes for \S \ref{42}}

The minimum given in (\ref{95}):

Clear[k, l, lambda, Q];

Q = Exp[-2*lambda]*(Sum[lambda^(2*k)/(k!)^2, \{k, 0, 250\}] +

     2*Sum[lambda^(2*k + l)/((k!)*(k + l)!), \{k, 0, 250\}, \{l, 1,
        IntegerPart[(2*lambda)^(1/2)]\}]);

QT = Table[Q, \{lambda, 0.5, 200, 1/2/10^3\}];

Min[QT]
\vskip 0.3cm
Output result: 0.564565

\subsection{Proof of Theorem \ref{pro2}}

For $\lambda>0$, define
$$
I_p^{'}(\lambda)=P\left\{|N_\lambda-E[N_\lambda]|< \sqrt{{\rm Var}(N_\lambda)}\right\}.
$$

\noindent (i) For $\lambda\in(0,\frac{3-\sqrt{5}}{2}\approx0.38197]$, we have $\lambda-\sqrt{\lambda}<0$ and $0<\lambda+\sqrt{\lambda}\le 1$. Then,
\begin{eqnarray*}
I_p^{'}(\lambda)=P\{N_\lambda=0\}=e^{-\lambda}\ge e^{-\frac{3-\sqrt{5}}{2}}\approx0.68252.
\end{eqnarray*}

\noindent (ii)  For $\lambda\in(\frac{3-\sqrt{5}}{2},1)$, we have $\lambda-\sqrt{\lambda}<0$ and $1<\lambda+\sqrt{\lambda}<2$. Then,
\begin{eqnarray*}
I_p^{'}(\lambda)=P\{N_\lambda=0\}+P\{N_\lambda=1\}=e^{-\lambda}(1+\lambda)>2e^{-1}\approx0.73576.
\end{eqnarray*}

\noindent (iii) For $\lambda=1$, we have
\begin{eqnarray*}
I_p^{'}(\lambda)=P\{N_\lambda=1\}=e^{-1}\approx0.36788.
\end{eqnarray*}

\noindent (iv) For $\lambda\in(1,\frac{7-\sqrt{13}}{2}\approx1.69722]$, we have $0<\lambda-\sqrt{\lambda}<1$ and $2<\lambda+\sqrt{\lambda}\le 3$. Then,
\begin{eqnarray*}
I_p^{'}(\lambda)=P\{N_\lambda=1\}+P\{N_\lambda=2\}=e^{-\lambda}\left(\lambda+\frac{\lambda^2}{2}\right).
\end{eqnarray*}

Let $f_1(\lambda)=e^{-\lambda}(\lambda+\frac{\lambda^2}{2})$. Since  $f_1'(\lambda)=e^{-\lambda}(1-\frac{\lambda^2}{2})$,  $f_1(\lambda)$ is increasing in $(1,\sqrt{2})$  and decreasing in $(\sqrt{2},\frac{7-\sqrt{13}}{2}]$. Note that $f_1(1)=1.5e^{-1}\approx0.55182$ and $f_1(\frac{7-\sqrt{13}}{2})\approx0.57476$. Then,
$$I_p^{'}(\lambda)>f_1(1)=1.5e^{-1}.$$

\noindent (v) For $\lambda\in(\frac{7-\sqrt{13}}{2},\frac{9-\sqrt{17}}{2}]$, we have $0<\lambda-\sqrt{\lambda}<1$ and $3<\lambda+\sqrt{\lambda}\le 4$. Then,
\begin{eqnarray*}
I_p^{'}(\lambda)=P\{N_\lambda=1\}+P\{N_\lambda=2\}+P\{N_\lambda=3\}=e^{-\lambda}\left(\lambda+\frac{\lambda^2}{2}+\frac{\lambda^3}{6}\right).
\end{eqnarray*}

Let $f_2(\lambda)=e^{-\lambda}(\lambda+\frac{\lambda^2}{2}+\frac{\lambda^3}{6})$. Since $f'_2(\lambda)=e^{-\lambda}(1-\frac{\lambda^3}{6})$,  $f_2(\lambda)$  is increasing in $(\frac{7-\sqrt{13}}{2},\sqrt[3]{6}\approx1.8171)$  and decreasing in $(\sqrt[3]{6},\frac{9-\sqrt{17}}{2}]$. Note that $f_2(\frac{9-\sqrt{17}}{2})\in(0.68335,0.68336)$ and $ f_2(\frac{7-\sqrt{13}}{2})\approx 0.72403$. Then,
$$I_p^{'}(\lambda)\ge f_2\left(\frac{9-\sqrt{17}}{2}\right)>0.68335.$$

\noindent (vi) For $\lambda\in(\frac{9-\sqrt{17}}{2},\frac{3+\sqrt{5}}{2}\approx2.618)$, we have $0<\lambda-\sqrt{\lambda}<1$ and $4<\lambda+\sqrt{\lambda}<5$. Then,
\begin{eqnarray*}
I_p^{'}(\lambda)=P\{N_\lambda=1\}+P\{N_\lambda=2\}+P\{N_\lambda=3\}+P\{N_\lambda=4\}=e^{-\lambda}\left(\lambda+\frac{\lambda^2}{2}+\frac{\lambda^3}{6}+\frac{\lambda^4}{24}\right).
\end{eqnarray*}

Let $f_3(\lambda)=e^{-\lambda}(\lambda+\frac{\lambda^2}{2}+\frac{\lambda^3}{6}+\frac{\lambda^4}{24})$. Since $f'_3(\lambda)=e^{-\lambda}(1-\frac{\lambda^4}{24})$,  $f_3(\lambda)$ is decreasing in $(\frac{9-\sqrt{17}}{2},\frac{3+\sqrt{5}}{2}\approx2.618)$. Note that $f_3(\frac{3+\sqrt{5}}{2})\in(0.80191,0.80192)$. Then,
$$I_p^{'}(\lambda)>0.80191.$$

\noindent (vii)  For $\lambda\in[\frac{3+\sqrt{5}}{2},\frac{11-\sqrt{21}}{2}\approx3.20871]$, we have $1\le\lambda-\sqrt{\lambda}<2$ and $4<\lambda+\sqrt{\lambda}\le5$. Then,
\begin{eqnarray*}
I_p^{'}(\lambda)=P\{N_\lambda=2\}+P\{N_\lambda=3\}+P\{N_\lambda=4\}=e^{-\lambda}\left(\frac{\lambda^2}{2}+\frac{\lambda^3}{6}+\frac{\lambda^4}{24}\right).
\end{eqnarray*}

Let $f_4(\lambda)=e^{-\lambda}(\frac{\lambda^2}{2}+\frac{\lambda^3}{6}+\frac{\lambda^4}{24})$. Since $f'_4(\lambda)=e^{-\lambda}(\lambda-\frac{\lambda^4}{24})$, $f_4(\lambda)$ is increasing in $[\frac{3+\sqrt{5}}{2},\sqrt[3]{24}\approx2.8845]$  and decreasing in $(\sqrt[3]{24},\frac{11-\sqrt{21}}{2}]$. Note that $f_4(\frac{3+\sqrt{5}}{2})\approx0.61094$ and $f_4(\frac{11-\sqrt{21}}{2})\in(0.60899,0.609)$. Then,
$$I_p^{'}(\lambda)>0.60899.$$

\noindent (viii) For $\lambda\in(\frac{11-\sqrt{21}}{2},4)$, we have $1<\lambda-\sqrt{\lambda}<2$ and $5<\lambda+\sqrt{\lambda}<6$. Then,
\begin{eqnarray*}
I_p^{'}(\lambda)&=&P\{N_\lambda=2\}+P\{N_\lambda=3\}+P\{N_\lambda=4\}+P\{N_\lambda=5\}\\
&=&e^{-\lambda}\left(\frac{\lambda^2}{2}+\frac{\lambda^3}{6}+\frac{\lambda^4}{24}+\frac{\lambda^5}{120}\right).
\end{eqnarray*}

Let $f_5(\lambda)=e^{-\lambda}(\frac{\lambda^2}{2}+\frac{\lambda^3}{6}+\frac{\lambda^4}{24}+\frac{\lambda^5}{120})$. Since $f'_5(\lambda)=e^{-\lambda}(\lambda-\frac{\lambda^5}{120})$,  $f_5(\lambda)$ is increasing in $(\frac{11-\sqrt{21}}{2},\sqrt[4]{120}\approx3.30975]$  and decreasing in $(\sqrt[4]{120},4)$. Note that $f_5(4)\in(0.69355,0.69356)$ and $f_5(\frac{11-\sqrt{21}}{2})\approx0.72353$. Then,
$$I_p^{'}(\lambda)>0.69355.$$

\noindent (ix) For $\lambda=4$, we have $\lambda-\sqrt{\lambda}=2$ and $\lambda+\sqrt{\lambda}=6$. Then,
\begin{eqnarray*}
I_p^{'}(\lambda)&=&P\{N_\lambda=3\}+P\{N_\lambda=4\}+P\{N_\lambda=5\}\\
&=&e^{-\lambda}\left(\frac{\lambda^3}{6}+\frac{\lambda^4}{24}+\frac{\lambda^5}{120}\right).
\end{eqnarray*}
Thus,
$$I_p^{'}(4)>0.547.$$

\noindent (x) For $\lambda\in(4, \frac{15-\sqrt{29}}{2}\approx4.80742]$, we have $2<\lambda-\sqrt{\lambda}<3$ and $6<\lambda+\sqrt{\lambda}\le7$. Then,
\begin{eqnarray*}
I_p^{'}(\lambda)&=&P\{N_\lambda=3\}+P\{N_\lambda=4\}+P\{N_\lambda=5\}+P\{N_\lambda=6\}\\
&=&e^{-\lambda}\left(\frac{\lambda^3}{6}+\frac{\lambda^4}{24}+\frac{\lambda^5}{120}+\frac{\lambda^6}{720}\right).
\end{eqnarray*}

Let $f_6(\lambda)=e^{-\lambda}(\frac{\lambda^3}{6}+\frac{\lambda^4}{24}+\frac{\lambda^5}{120}+\frac{\lambda^6}{720})$. Since $f'_6(\lambda)=e^{-\lambda}(\frac{\lambda^2}{2}-\frac{\lambda^6}{720})$,  $f_6(\lambda)$ is decreasing in $(4, \frac{15-\sqrt{29}}{2}]$ and $f_6(\frac{15-\sqrt{29}}{2})\in(0.64792,0.64793)$. Then,
$$I_p^{'}(\lambda)>0.64792.$$

\noindent (xi) For $\lambda> \frac{15-\sqrt{29}}{2}$, let $g_1(\lambda):=P\{N_\lambda< \lambda+\sqrt{\lambda}\}$ and $g_2(\lambda):=P\{N_\lambda\le \lambda-\sqrt{\lambda}\}$. Then,
$$I_p^{'}(\lambda)=g_1(\lambda)-g_2(\lambda).$$

\noindent (xia) Suppose $\lambda=n$ for some $n\in\mathbb{N}$. Let $Y_1,Y_2,\dots,Y_n$ be i.i.d. Poisson(1) random variables. We have $N_n=Y_1+Y_2+\cdots+Y_n\ \ {\rm in\ distribution}$.
Note that $C\rho<0.7656(1+2e^{-1})<1.328898$. Then,  by (\ref{K1}), we obtain that for any $ n\ge108$,
 $$P\left\{N_n\le n+\sqrt{n}\right\}\ge \Phi(1)-\frac{1.328898}{\sqrt{108}}>0.7134.$$

\noindent (xib) For $n\ge108$, consider $\lambda\in (n,n+1)$. We have
 $$P\{N_\lambda< \lambda+\sqrt{\lambda}\}\ge P\{N_\lambda\le n+\sqrt{n}\}> P\{N_{n+1}\le n+\sqrt{n}\}.$$
Then, by (\ref{K2}), we get
 \begin{eqnarray*}
P\left\{N_{n+1}\le n+\sqrt{n}\right\}&\ge& \Phi\left(\frac{\sqrt{n}-1}{\sqrt{n+1}}\right)-\frac{C\rho}{\sqrt{n+1}}\\
&\ge & \Phi\left(\frac{\sqrt{108}-1}{\sqrt{108+1}}\right)-\frac{C\rho}{\sqrt{108+1}}\\
 &>&\Phi(0.8996)-0.1273\\
&>&0.6886.
 \end{eqnarray*}
Thus,
\begin{equation}\label{g111}
 g_1(\lambda)>0.6886,\ \ \ \ \forall \lambda\ge108.
\end{equation}

\noindent (xic)  Suppose $\lambda=n$ for some $n\in\mathbb{N}$. For any $n\ge108$, we have  $\Phi(-1)+\frac{C\rho}{\sqrt{108}}< 0.2866$. Then, by (\ref{K3}), we get
 $$P\left\{N_n\le n-\sqrt{n}\right\}<0.2866.$$

\noindent (xid) Consider $\lambda\in (n,n+1)$ for $n\ge108$.  We have
 $$P\{N_\lambda\le \lambda-\sqrt{\lambda}\}< P\{N_n\le \lambda-\sqrt{\lambda}\}\le P\{N_{n}\le n+1-\sqrt{n+1}\}.$$
Then, by (\ref{K4}), we get
\begin{eqnarray*}
 P\{N_{n}\le n+1-\sqrt{n+1}\}&\le&  \Phi\left(\frac{1-\sqrt{108+1}}{\sqrt{108}}\right)+\frac{C\rho}{\sqrt{108}}\\
&<&\Phi(-0.9083)+0.1279\\
&<&0.3093.
\end{eqnarray*}
 Thus,
\begin{equation}\label{g222}
  g_2(\lambda)<0.3093, \ \ \ \ \forall \lambda\ge108.
\end{equation}
Therefore, by (\ref{g111}) and (\ref{g222}),  we find that for $\lambda\ge108$,
 $$I_p^{'}(\lambda)=g_1(\lambda)-g_2(\lambda)>0.6886-0.3093=0.3793>e^{-1}.$$

\noindent (xie) Finally, we consider the case that  $\lambda\in(\frac{15-\sqrt{29}}{2},108)$.

\noindent (xie1) For $n\in\mathbb{N}$, let $\lambda_{n,1}^{(1)},\lambda_{n,2}^{(1)}$ be the solutions of $\lambda+\sqrt{\lambda}=n$, where
$$\lambda_{n,1}^{(1)}=\frac{2n+1-\sqrt{4n+1}}{2},\ \ \ \ \lambda_{n,2}^{(1)}=\frac{2n+1+\sqrt{4n+1}}{2}.$$
Since $\lambda<n$, the equation has a unique solution, which is denoted by $ \lambda_{n}^{(1)}=\frac{2n+1-\sqrt{4n+1}}{2}<n-1$.

For fixed $x$, $P\{N_\lambda< x\}$ is decreasing with respect to $\lambda$. Then,  when $\lambda+\sqrt{\lambda}\in(n,n+1]$, $P\{N_\lambda<\lambda+\sqrt{\lambda}\}=P\{N_\lambda\le n\}$  is decreasing for $\lambda\in (\lambda_{n}^{(1)},\lambda_{n+1}^{(1)}]$. Thus,  for $\lambda\in(\lambda_{n}^{(1)},\lambda_{n+1}^{(1)}]$, we have
$$g_1(\lambda)\ge P\{N_{\lambda_{n+1}^{(1)}}\le n\}.$$
Note that $\frac{15-\sqrt{29}}{2}=\lambda_{7}^{(1)}$ and $ 108\in(\lambda_{118}^{(1)},\lambda_{119}^{(1)})$. For $\lambda\in(\frac{15-\sqrt{29}}{2},108)$,  we use Matlab to find the minimal value of $P\{N_{\lambda_{n}^{(1)}}\le n-1\}$ for $ n\in[8,119]$, which is $P\{N_{\lambda_{8}^{(1)}}\le 7\}\in(0.79345,0.79346)$. Hence,  we have
\begin{equation}\label{g13}
  g_1(\lambda)>0.79345,\ \ \ \ \lambda\in\left(\frac{15-\sqrt{29}}{2},108\right).
\end{equation}

\noindent (xie2) For $n\in\mathbb{N}$, let $\lambda_{n,1}^{(2)},\lambda_{n,2}^{(2)}$ be the solutions of $\lambda-\sqrt{\lambda}=n$, where
$$\lambda_{n,1}^{(2)}=\frac{2n+1-\sqrt{4n+1}}{2},\ \ \ \ \lambda_{n,2}^{(2)}=\frac{2n+1+\sqrt{4n+1}}{2}.$$
Since $\lambda>n$, the equation has a unique solution, which is denoted by $\lambda_{n}^{(2)}=\frac{2n+1+\sqrt{4n+1}}{2}>n+1$.

For fixed $x$, $P\{N_\lambda\le x\}$ is decreasing with respect to $\lambda$. Then, when $\lambda-\sqrt{\lambda}\in[n,n+1)$, $P\{N_\lambda\le\lambda-\sqrt{\lambda}\}=P\{N_\lambda\le n\}$ is decreasing for $\lambda\in [\lambda_{n}^{(2)},\lambda_{n+1}^{(2)})$. Thus, for $\lambda\in [\lambda_{n}^{(2)},\lambda_{n+1}^{(2)})$, we have
$$g_2(\lambda)\le P\{N_{\lambda_{n}^{(2)}}\le n\}.$$
Note that $\frac{15-\sqrt{29}}{2}\in(\lambda_{2}^{(2)},\lambda_{3}^{(2)})$ and $108\in(\lambda_{97}^{(2)},\lambda_{98}^{(2)})$. For $\lambda\in(\frac{15-\sqrt{29}}{2},108)$, we consider the maximal value of $g_2(\frac{15-\sqrt{29}}{2})$ and  $P\{N_{\lambda_{n}^{(2)}}\le n\}$ for $ n\in[3,97]$. We have $g_2(\frac{15-\sqrt{29}}{2})=P\{N_{\frac{15-\sqrt{29}}{2}}\le2\}\approx0.14184$. By virtue of Matlab, we find the maximal value of  $P\{N_{\lambda_{n}^{(2)}}\le n\}$ for $ n\in[3,97]$, which is $P\{N_{\lambda_{3}^{(2)}}\le 3\}\in(0.22506,0.22507)$.  Hence, we have
\begin{equation}\label{g23}
  g_2(\lambda)<0.22507,\ \ \ \ \lambda\in\left(\frac{15-\sqrt{29}}{2},108\right).
\end{equation}

By (\ref{g13}) and (\ref{g23}), we conclude that for $\lambda\in(\frac{15-\sqrt{29}}{2},108)$,
$$
  I_p^{'}(\lambda)=g_1(\lambda)-g_2(\lambda)>0.79345-0.22507=0.56838>e^{-1}.
$$
The proof of Theorem \ref{pro2} is therefore complete based on all the above cases. \hfill \fbox

\vskip 0.5cm
{ \noindent {\bf\large Acknowledgements}\quad This work was supported by the National Natural Science Foundation of China (Nos. 12161029 and 12171335), the National Natural Science Foundation of Hainan Province (No. 121RC149), the Science Development Project of Sichuan University (No. 2020SCUNL201)  and the Natural Sciences and Engineering Research Council of Canada (No. 4394-2018).

\end{document}